\documentclass[12pt]{amsart}
\usepackage{latexsym, amssymb, amsmath}
 \usepackage{eucal}

    \newtheorem{rema}{Remark}[section]
    \newtheorem{propo}[rema]{Proposition}

   \newtheorem{theo}[rema]{Theorem}
   \newtheorem{def-theo}[rema]{Definition-Theorem}
 \newtheorem{conj}[rema]{Conjecture}
   
    \newtheorem{lemma}[rema]{Lemma}

  \newtheorem{rmk}[rema]{Remark}
\newtheorem{quest}[rema]{Question}
      
	\newcommand{\nno}{\nonumber}

	\newcommand{\p}{\partial}

 \newcommand{\pf}{{\it Proof:}\hspace{2ex}}
 
 \newcommand{\epfv}{\hspace{1em}$\Box$\vspace{1em}}

\newcommand{\bC}{{\mathbb C}}
\newcommand{\bZ}{{\mathbb Z}}

\newcommand{\bQ}{{\mathbb Q}}
\newcommand{\bN}{{\mathbb N}}
\newcommand{\bT}{{\mathbb T}}

\newcommand{\BQ}{\begin{eqnarray}}
\newcommand{\EQ}{\end{eqnarray}}
\newcommand{\BQn}{\begin{eqnarray*}}
\newcommand{\EQn}{\end{eqnarray*}}
\newcommand{\BL}{\begin{align}}
\newcommand{\EL}{\end{align}}
\newcommand{\BLn}{\begin{align*}}
\newcommand{\ELn}{\end{align*}}
\newcommand{\BA}{\begin{align}}
\newcommand{\EA}{\end{align}}
\newcommand{\BAn}{\begin{align*}}
\newcommand{\EAn}{\end{align*}}

\newcommand{\cP}{\mathcal P}

\title[ Recurrent Inversion Formulas]
{Recurrent Inversion Formulas}

    \author{Wenhua Zhao}      
    \begin{document}

\begin{abstract}
Let $F(z)=z-H(z)$ with $o(H(z))\geq 2$ be a formal map from $\bC^n$ to $\bC^n$
and $G(z)$ the formal inverse of $F(z)$.
In this paper, we fist study the deformation 
$F_t(z)=z-tH(z)$ and its formal 
inverse map $G_t(z)$. We then derive 
two recurrent 
formulas for the formal inverse $G(z)$. The first formula 
in certain situations 
provides a more efficient method for the calculation of $G(z)$ 
than other well known inversion formulas.
The second one is differential free but 
only works when $H(z)$ is 
homogeneous of degree $d\geq 2$. Finally, we reveal 
a close relationship 
of the inversion problem with a Cauchy problem of a PDE. 
When the Jacobian matrix $JF(z)$ is symmetric, 
the PDE coincides with the $n$-dimensional inviscid 
Burgers' equation in Diffusion theory.
\end{abstract}

   

 \bibliographystyle{alpha}
    \maketitle
\tableofcontents

\renewcommand{\theequation}{\thesection.\arabic{equation}}
\renewcommand{\therema}{\thesection.\arabic{rema}}
\setcounter{equation}{0}
\setcounter{rema}{0}
\setcounter{section}{0}

\section{\bf Introduction}\label{S1}

Let $F(z)=z-H(z)$  be a formal map from $\bC^n$ to 
$\bC^n$ with $o(H(z))\geq 2$
and $G(z)$ the formal inverse of $F(z)$, i.e. 
$z=F(G(z))=G(F(z))$. The formulas which directly or indirectly 
give the formal inverse $G(z)$ 
are called inversion formulas in literature. 
There have been many different versions
of inversion formulas.
The first inversion formula in history 
was the Lagrange's inversion formula given 
by L. Lagrange \cite{L} in 1770, 
which provides a formula 
to calculate all coefficients of $G(z)$ 
for the one-variable case. This formula was generalized to multi-variable
cases by I. G. Good \cite{Go} in 1965. Jacobi \cite{J1} in 1830 also 
gave an inversion formula for the cases $n\leq 3$ and later \cite{J2} 
for the general case. This formula now is 
called the Jacobi's inversion formula. 
Another inversion formula is Abhyankar-Gurjar inversion formula, 
which was first proved by Gurjar in 1974 (unpublished) 
and later Abhyankar  \cite{A} gave a simplified proof.
By using Abhyankar-Gurjar inversion formula, 
H. Bass, E. Connell and D. Wright \cite{BCW} and D. Wright \cite{W2} 
proved the so-called Bass-Connell-Wright's tree expansion formula. 
Recently, in \cite{WZ}, this formula has been generalized
to a tree expansion formula for the formal flow $F_t(z)$ 
of $F(z)$ which provides a uniform  
formula for all the powers $F^{[m]}(z)$ $(m\in \bZ)$ of $F(z)$.  
There are also many other inversion formulas in literature, 
see \cite{Ge}, \cite{W3} and references there.

One of the motivations of seeking inversion formulas comes from their 
important applications in enumerative combinatorics of rooted trees. 
See, for example, \cite{St1}, \cite{Ge} and references there.
Another motivation comes from the study of 
the well known Jacobian conjecture.
Recall that the Jacobian conjecture claims that,  
any polynomial map $F(z)$ from $\bC^n$ to $\bC^n$
with the Jacobian $j(F)=\det (\frac {\p F_i}{\p z_j})=1$ 
is an automorphism 
of $\bC^n$ and its inverse $G(z)$ is also a polynomial map. 
The Jacobian conjecture was first proposed by Keller \cite{K} 
in 1939. It is now still 
open even for the two variable case. 
For the history and some well known results 
of the Jacobian conjecture, 
see the classical paper \cite{BCW}, \cite{E} 
and references there. 

In this paper, we fist study the deformation 
$F_t(z)=z-tH(z)$ and its formal 
inverse map $G_t(z)$. We then derive 
two recurrent 
formulas for the formal inverse $G(z)$. 
The first formula 
in certain situations 
provides a more efficient method for the calculation of $G(z)$ 
than other well known inversion formulas.
The second one is differential free but 
only works when $H(z)$ is 
homogeneous of degree $d\geq 2$. Finally, we reveal 
a close relationship 
of the inversion problem with a Cauchy problem of a PDE. 
When the Jacobian matrix $JF(z)$ is symmetric, 
the PDE coincides with the $n$-dimensional inviscid
Burgers' equation in Diffusion theory.

The arrangement of this paper is as follows.
In Section \ref{S2}, we recall several well known inversion 
formulas in literature and derive the formulas for the formal inverse  
$G(z)$ if they are not given directly. 
In Section \ref{S3}, we mainly study the deformation 
$F_t(z)=z-tH(z)$ and its formal 
inverse map $G_t(z)$. The main results are 
Theorem \ref{T3.2} and Proposition \ref{P3.3}.
Theorem \ref{T3.2} reveals a close relationship 
between the formal inverse $G_t(z)$ and a Cauchy problem 
of PDE which has a similar form as 
the $n$-dimensional inviscid Burgers' equation. 
Proposition \ref{P3.3} give us 
the first recurrent inversion formula.
In this section, we also derive some 
interesting consequence from the main results above.
One of theorem is Proposition \ref{New-P} which 
claims that the formal inverse $G(z)=z+H(z)$ 
if and only if $JH\cdot H=0$. As an immediate consequence, 
we derive the well known theorem of H. Bass, E. Connell and 
D. Wright \cite{BCW}, which says that the Jacobian conjecture 
is true when $H(z)$ is homogeneous and $J(H)^2=0$. 
One purpose that we include the proof for this theorem 
is to explore why it is much more difficult to prove 
the Jacobian conjecture under the same conditions as 
above except we have $J(H)^k=0$ for some $k\geq 3$ 
instead of $J(H)^2=0$.  
Finally, in  Section \ref{S4}, 
we derive our second recurrent inversion formula.
This formula is a generalization of the recurrent formula 
proved by Dru\.zkowski 
\cite{D2} for the case $d=3$.
Our approach also gives 
the involved multi-linear form explicitly.

Finally, two remarks are as follows. 
First, we will fix $\bC$ as our base field.
But all results, formulas as well as their proofs
given in this paper hold or work equally well 
for formal power series over any $\bQ$-algebra. 
Secondly, for convenience,
we will mainly work on the setting of 
formal power series over $\bC$. But, 
for polynomial maps or local analytic maps, 
all formal maps or 
power series involved 
in this paper are also
locally convergent.
This can be easily seen either from 
the fact that any local analytic map with non-zero Jacobian 
at the origin has a locally convergent inverse, or from 
the well-known Cauchy-Kowaleskaya theorem (See \cite{R}, for example.) 
in PDE.

The author would like to thank Professor David Wright 
for informing the author
that the recurrent formulas presented here are new in the literature 
and also for encouraging him to write this note. The author 
would also like to thank Professor Quo-Shin Chi 
for discussions on some PDE's involved in this paper. 

\renewcommand{\theequation}{\thesection.\arabic{equation}}
\renewcommand{\therema}{\thesection.\arabic{rema}}
\setcounter{equation}{0}
\setcounter{rema}{0}

\section{\bf Inversion Formulas}\label{S2}

In this section, we review the inversion formulas of Lagrange, 
Jacobi, Abhyankar-Gurjar, Bass-Connell-Wright and also the formula 
developed in \cite{WZ} for the formal flow $F_t(z)$ of $F(z)$, 
which encodes formulas for all 
powers $F^{[m]}(z)$ $(m\in \bZ)$. There are also many 
other versions of inversion formulas in literature. They 
are more or less in the same favor as some of inversion formulas above. 
We refer readers to \cite{Ge}, \cite{W3} and 
references for other inversion formulas.

First we fix the following notation.

We let $z=(z_1, z_2, \cdots, z_n)$ and $\bC[z]$ (resp. $\bC[[z]]$) 
the algebra of polynomials (resp. formal power series) in $z$. For any 
${\bf k}\in \bZ^n$ and Laurent series $h(z)$, we denote 
by $[z^{\bf k}]h(z)$ the coefficient of $[z^{\bf k}]$ in $h(z)$. 
In particular, we set 
\BQ
\text{Res}_z h(z)=[z_1^{-1}z_2^{-1}\cdots z_n^{-1}] h(z)
\EQ

In this paper, we always denoted by $F(z)=(F_1(z), \cdots, F_n(z))$ 
a formal map from $\bC^n$ to itself 
with the form $F(z)=z-H(z)$ and $o(H(z))\geq 2$. 
The notation $G(z)$ always denotes the formal inverse of $F(z)$
and $j(F)$ the Jacobian $\det (\frac {\p F_i}{\p z_j})=1$ of $F(z)$.

We start with the Lagrange's multi-variable inversion formula. The version 
of the formula we quote here is given by Good in \cite{Go}.

\begin{theo} $(\text{\bf Lagrange's Inversion Formula})$ 

Let $f(z)=(f_1(z), f_2(z), \cdots, f_n(z))\in \bC[[z]]^{\times n}$ with non-zero
constant terms. Then there exist a unique 
$g(z)=(g_1(z), \cdots, g_n(z))\in \bC[[z]]^{\times n}$ such that
\BQ\label{E2.2}
g_i(z)=z_i f_i(g(z))
\EQ
for any $i=1, 2, \cdots n$.
Furthermore, for any formal Laurent series $\phi (z)$ and 
${\bf k}\in \bZ^n$, 
we have 
\begin{align} 
[z^{{\bf k}}]\frac {\phi (g)}{\det (\delta_{i, j}-z_i
\frac {\p f_i}{\p z_j}(g) )} =[w^{{\bf k}}] \phi(w)f^{{\bf k}}(w) \\
[z^{{\bf k}}]\phi(g) = [w^{{\bf k}}]\det (\delta_{i, j}-\frac {w_i}{f_i(w)}
\frac {\p f_i}{\p z_j}(w) ) \phi(w)f^{{\bf k}}(w)\label{L-Inv}
\end{align}
\end{theo}

Let us see how to use the formulas above to calculate the formal inverse $G(z)$
of $F(z)$. To do this, we assume that 
\BQ\label{E2.5}
H_i(z)=z_i h_i(z) \quad (1\leq i\leq n)
\EQ
for some $h_i(z)\in \bC[[z]]$. 

We choose $f_i(z)=\frac 1{1-h_i(z)}$. Hence,
$\frac {z_i}{f_i(z)}=F_i(z)$ and by
Eq.\,(\ref{E2.2}), we have
\BQn
z_i=\frac {z_i}{f_i} (g(z))=F_i(g(z))
\EQn
i.e. $g(z)$ is the formal inverse of $F(z)$, therefore $g(z)=G(z)$.

To calculate $G(z)$, hence,
it's enough to calculate $[z^{\bf k}]G_i(z)$ for any 
$1\leq i\leq n$ and ${\bf k}\in \bN^n$. For any fixed $i$, we choose
$\phi (z)=z_i$. By Eq.\,(\ref{L-Inv}), we get
\BQ
[z^{{\bf k}}] G_i(z) = [w^{{\bf k}}]\det (\delta_{i, j}-\frac {w_i}{f_i(w)}
\frac {\p f_i}{\p z_j}(w) ) w_i f^{{\bf k}}(w)
\EQ

For the case when $H_i(z)$ is not of the form (\ref{E2.5}).
The Lagrange's inversion formula does not provide a direct method for
the calculation of $G(z)$. But, one can derive from 
the Lagrange's inversion formula
the following two inversion formulas, 
which will provide formulas for $G(z)$ 
in the general case. (See, for example, \cite{Ge}).

The next inversion formula was first proved by Jacobi \cite{J1} for 
$n\leq 3$ and later in \cite{J2} for the general case.

\begin{theo} $(\text{\bf Jacobi's Inversion Formula})$ 

Let $f(z)=(f_1(z), f_2(z), \cdots, f_n(z))$ be a sequence of formal
series in $z$. Let $\phi(z)$ be any Laurent series. Then
\BQ\label{J-Inv}
\mbox{Res}_w\phi(w)=\mbox{Res}_z j(f) \phi(f(z))
\EQ
\end{theo}

To get the formal inverse $G(z)$ of $F(z)$ by using the 
Jacobi's inversion formula, 
we can choose $f(z)=F(z)$. For each $1\leq i\leq n$ and ${\bf k}\in \bN^n$, we
choose $\phi (z)=z^{-{\bf k}-{\bf 1}} G_i(z)$, 
where ${\bf 1}=(1, 1, \cdots, 1)$. 
Then Formula (\ref{J-Inv}) gives us
\BQ
[z^{\bf k}]G_i(z)= \text{Res}_w w^{-{\bf k}-{\bf 1}}G_i (w)=
\mbox{Res}_z j(F) F^{-{\bf k}-{\bf 1}}(z) z_i 
\EQ

Hence, by changing ${\bf k}\in \bN^n$, 
we can calculate $G(z)$ completely. 

The first direct inversion formula was proved by Gurjar (unpublished) 
and a simplified proof was later given by Abhyankar \cite{A}.

\begin{theo} $(\text{\bf Abhyankar-Gurjar's Inversion Formula})$

Let $F(z)=z-H(z)$ with $o(H)\geq 2$. Then for any Laurent series $\phi(z)$, 
we have
\BQ\label{A-G}
\phi (G)=\sum_{{\bf m}\in \bN^n} \frac { D^{\bf m} } {{\bf m}!} 
\phi(z) j(F) H^{\bf m} 
\EQ
where $D^{\bf m}=D_1^{m_1}D_2^{m_2}\cdots D_n^{m_n}$ for any 
${\bf m}=(m_1, \cdots, m_2)\in \bN^n$ and $D_i=\frac {\p}{\p z_i}$ for any
$1\leq i\leq n$.
\end{theo}

Note that, if we choose $\phi(z)=z_i$ $(1\leq i\leq n)$, 
we get $G_i(z)$ by Formula (\ref{A-G}).

By using Abhyankar-Gurjar's inversion formula,  H. Bass, E. Connell and D. Wright 
\cite{BCW} proved the so-called Bass-Connell-Wright's tree expansion inversion 
formula for the case when $H(z)$ is homogeneous and later D. Wright \cite{W2}
proved that the same formula also holds in the general case. A 
totally different proof of this formula was also given in \cite{WZ}.

In order to recall Bass-Connell-Wright's inversion 
formula and the formula for the formal flow $F_t(z)$ of $F(z)$ 
developed in \cite{WZ}, we need the following notation. 

By a {\it rooted tree} $T$ we mean a finite
connected and simply connected 
graph with one vertex designated as its {\it root}, denoted by $rt_T$.
In a rooted tree
there are natural ancestral relations between vertices.  We say a
vertex $w$ is a {\it child} of vertex $v$ if the two are connected by an
edge and $w$ lies further from the root than $v$. We denote by $v^+$ 
the set of all its children.
When we speak of isomorphisms between rooted trees, we will
always mean root-preserving isomorphisms. 
We denote by $\bT$ (resp. $\bT_m$) the set of equivalent classes of
rooted trees (resp. rooted trees with $m$ vertices).

For $T\in\bT$, we denote by $V(T)$ the set of vertices of $T$ and
$|T|=|V(T)|$.
A {\it labeling} of $T$ in the set $\{1,\ldots,n\}$ is
a function $l:V(T)\to\{1,\ldots,n\}$. A rooted tree $T$ with a
labeling $l$ is called a {\it labeled rooted tree,} denoted $(T, l)$.
Given such, and given $F=z-H$ as above, we make the following
definitions, for $v\in V(T)$:
\begin{enumerate}
\item $H_v=H_{l(v)}$.
\item $D_v=D_{l(v)}$.
\item $D_{v^+}=\prod_{w\in v^+} D_w$.
\item $P_{T,l}=\prod_{v\in V(T)} D_{v^+}H_v$.
\end{enumerate}
Finally, we define systems of power series
 $\cP_T=(\cP_{T,1},\ldots, \cP_{T,n})$
by setting
\BQ
\cP_{T,i}=\frac{1}{|\text{Auto(T)}|} 
\sum_{\substack{{l:V(T)\to\{1,\ldots,n\}}\\f(\text{rt}_T)=i}}
P_{T,l}
\EQ
for $i=1,\ldots,n$, where the sum above runs over all
labelings of $T$ having a fixed label for the root.

\begin{theo} $(\text{\bf Bass-Connell-Wright's Inversion Formula})$ 

With the fixed notation above, the formal inverse $G(z)$ of $F(z)$ is given by
\BQ\label{BCW-Inv}
G(z)=z+\sum_{T\in \bT} {\mathcal P}_T(z) 
\EQ
\end{theo}

For any $m\in \bZ$, we define the $m^{th}$-power $F^{[m]}(z)$ of $F(z)$ by
\begin{align}
F^{[m]}(z)&=\underbrace {F\circ F\circ \cdots \circ F}_{m \, \,  \text{copies} }(z)
\text{\quad if $m\geq 0$;}\label{F-m>0} \\
F^{[m]}(z)&=G^{[-m]}(z)\text{\quad if $m< 0$.}\label{F-m<0}
\end{align}

Considering all the efforts deriving formulas for the formal inverse $G(z)$ of 
the formal map $F(z)$, one may ask if there are some uniform formulas for 
all the powers $F^{[m]}(z)$ $(m\in \bZ)$ of $F(z)$. 
This question was answered in \cite{WZ} by deriving a formula 
for the formal flow $F_t(z)$, where $F_t(z)$ is 
the unique 1-parameter subgroup with $F_{t=1}(z)=F(z)$
of the formal automorphisms 
of $\bC^n$. The formula is derived by the D-log formulation \cite{Z1} 
of $F_t(z)$ and similar technics in \cite{BCW}.

\begin{theo} \cite{WZ} 
There is a unique sequence $\{\Psi_T(t) |T\in \bT\}$ such that
\BQ \label{W-Z}
F_t(z)=z+\sum_{T\in \bT} (-1)^{|T|} \Psi_T (t) {\mathcal P}_T(z) 
\EQ
\end{theo}

For some properties and a computational algorithm for 
the polynomials $\Psi_T(t)$, see \cite{WZ}.
Based on certain properties of $\Psi_T(t)$ proved in \cite{WZ}, 
J. Shareshian \cite{Sh} pointed out to us that
the  polynomials $\Psi_T(t)$ coincide with strict order polynomials 
$\bar{\Omega}(T, m)$ of rooted trees, 
which we will explain briefly below.

Note that, for any rooted tree $T$, the graph structure induces 
a natural partial order on the set $V(T)$ of vertices of $T$ with the root 
$rt_T$ of $T$ 
serving as the unique minimum element. Hence, with this partial order, 
the set $V(T)$ becomes a poset (partially ordered set). We will 
still use the same notation
$T$ to denote this poset. For any $m\geq 1$, we denote by $[m]$ the totally  
ordered set $\{1, 2, \cdots, m\}$. We say a map $\sigma$ from a finite poset
$P$ to $[m]$ is {\it strict order preserving} if 
$\sigma (x) > \sigma (y)$ in $[m]$ implies $x> y$ in $P$.
The following definition and theorem is well known in enumerative combinatorics. 
See, for example, \cite{St1}

\begin{def-theo}
For any finite poset $P$, there exists a unique polynomial $\bar{\Omega}(T, t)$
such that, for any $m\geq 1$, $\bar{\Omega}(T, m)$ equals to the number of 
order-preserving maps from $P$ to the totally ordered set $[m]$.
\end{def-theo}

\begin{theo} \cite{Sh}
For any rooted tree $T$, we have
\BQ
\Psi_T(t)=\bar \Omega (T, t)
\EQ
\end{theo}

\begin{rmk}
$(a)$ Note that the formula $(\ref{W-Z})$ provides a uniform formula 
for all powers $F^{[m]}(z)$ $(m\in \bZ)$ by setting $t=m$.

$(b)$ It is well known (See \cite{St1}.) that 
$\bar\Omega (T, -1)=(-1)^{|T|}$.  
For a direct proof of this fact, see \cite{WZ}. Hence,
by setting $m=-1$, 
the formula (\ref{W-Z}) becomes
Bass-Connell-Wright's inversion formula 
$(\ref{BCW-Inv})$.
\end{rmk}

\renewcommand{\theequation}{\thesection.\arabic{equation}}
\renewcommand{\therema}{\thesection.\arabic{rema}}
\setcounter{equation}{0}
\setcounter{rema}{0}

\section{\bf The First Recurrent Inversion Formula}\label{S3}

In this section, we first study a deformation 
of formal maps from which we then
derive our first recurrent inversion formula. 
Comparing with other well-known inversion formulas,
the recurrent inversion formula 
in certain situations provides 
a more efficient method 
 for the calculation of 
formal inverse maps. 
We also discuss a close relationship 
between the inversion problem and 
a Cauchy problem of a PDE, 
see Eq.\,(\ref{PDE}) and (\ref{PDE-B}). 

We start with the deformation $F_t(z)=z-tH(z)$ and 
let $G_t(z)$ be its formal inverse. 
Note that, we can always 
write the formal 
inverse $G_t(z)$ of $F_t(z)$ as 
$G_t(z)=z+tN_t(z)$ for a unique 
$N_t(z)\in \bC [t][[z]]^{\times n}$ 
and $\mbox{o} (N_t(z))\geq 2$.
This can be easily proved by using
any well-known reversion formulas, 
for example, Abhyankar-Gurjar inversion formula \cite{A} or
the Bass-Connell-Wright tree expansion formula \cite{BCW}.

\begin{lemma}\label{L3.1}
Let $F(z)$, $F_t(z)$, $G_t(z)$ and $N_t(z)$ be as above. Then we have
\BQ
N_t(F_t(z))&=&H(z), \label{E3.1}\\
H(G_t)&=& N_t(z), \label{E3.2}\\
JN_t(F_t)&=&JH(I-tJH)^{-1}=\sum_{k=1}^\infty JH^k(z)t^{k-1}.\label{MainEQ}
\EQ
In particular, for any $m\geq 1$, $JH^m(z)=0$ if and only if $JN_t^m(z)=0$.
\end{lemma}
\pf Since $z=G_t(F_t)$, we have
\BQn
 z &=& F_t(z)+tN_t(F_t(z)) \\
  z &=& z-tH(z)+tN_t(F_t(z))
\EQn
Therefore, 
\BQn
H(z)=N_t(F_t(z)),
\EQn
which is Eq.\,(\ref{E3.1}). Composing the both sides of Eq.\,(\ref{E3.1}) 
with $G_t(z)$, 
we get Eq.\,(\ref{E3.2}).

Now we show Eq.\,(\ref{MainEQ}). 
First,  by the fact $JG_t(F_t(z))=JF_t^{-1}(z)$, we have
\begin{align*}
I+tJN_t(F_t)&=(I-tJH)^{-1},\nno \\
tJN_t(F_t)= I-(I-&tJH)^{-1}=tJH(I-tJH)^{-1},\nno \\
JN_t(F_t)=JH(I-&tJH)^{-1}=\sum_{k=1}^\infty JH^k(z)t^{k-1}.
\end{align*}
Hence, we have Eq.\,(\ref{MainEQ}). 

Now, for any $m\geq 1$, by Eq.\,(\ref{MainEQ}), we have
\begin{align}
JN_t^m (F_t)=JH^m(I-&tJH)^{-m}.
\end{align}
Since $(I-tJH)^{-m}$ is invertible and $F_t(z)$ 
is an automorphism of $\bC[t][[z]]$, 
we have $JH^m(z)=0$ if and only if $JN_t^m(z)=0$.
\epfv

\begin{theo}\label{T3.2}
Let
$N_t(z)\in \bC [t][[z]]^{\times n}$ with
$\mbox{o} (N_t(z))\geq 2$. Then $G_t(z)=t+tN_t(z)$ is the formal inverse of
$F_t(z)=z-tH(z)$ if and only if $N_t(z)$ is the unique solution of 
the Cauchy problem of the following partial differential equation
\BQ
 &{}& \frac {\p N_t}{\p t}=JN_tN_t \label{PDE}\\
 &{}& N_{t=0}(z)=H(z) \label{PDE-B}
\EQ
where $JN_t$ is the Jacobian matrix of $N_t(z)$ with respect to $z$. 
\end{theo}
\pf
By applying $\frac {\p}{\p t}$ to the both sides of Eq.\,(\ref{E3.1}), we get
\BQn
0&=&\frac{\p N_t (F_t)}{\p t}\\
&=&\frac{\p N_t}{\p t}(F_t)+
JN_t (F_t)\frac{\p F_t}{\p t}\\
&=& \frac{\p N_t}{\p t}(F_t)-
JN_t (F_t)H 
\EQn
Therefore,
\BQn
 \frac{\p N_t}{\p t}(F_t)=
JN_t (F_t)H
\EQn
Composing with $G_t$, we get
\BQn
\frac{\p N_t}{\p t}=
JN_tH(G_t)=JN_t N_t
\EQn

Note that $G_{t=0}(z)=z$, for it is the formal inverse of 
$F_{t=0}(z)=z$.  Eq.\,(\ref{PDE-B}) 
follows immediately from Eq.\,(\ref{E3.2}) by setting $t=0$.

Note that the power series solution with respect to $t$ and $z$ 
of Eq.\,(\ref{PDE}) and (\ref{PDE-B}) 
is unique, hence, conversely, the
theorem is also true. 
\epfv

We define the sequence $\{ N_{[m]}(z) | m\geq 0\}$ by setting
$N_{[0]}(z)=z$ and writing 
$N_t(z)= \sum_{m=1}^\infty t^{m-1} N_{[m]}(z)$.

\begin{propo} \label{P3.3}
Let $N_t(z)= \sum_{m=1}^\infty t^{m-1} N_{[m]}(z)$
be the unique solution of  Eq. $(\ref{PDE})$ and $(\ref{PDE-B})$. Then
\BQ
N_{[1]} &=& H  \label{N1}\\
N_{[m]} &=& \frac 1{m-1} \sum_{\substack{k+l=m\\ k, l\geq 1}} 
JN_{[k]}\cdot N_{[l]}  \label{Nm}
\EQ
for any $m\geq 2$.
\end{propo}

\pf First, 
Eq.\,(\ref{N1}) follows immediately from Eq. $(\ref{PDE-B})$. 
Secondly, by Eq.\,(\ref{PDE}), we have

\BQn
\sum_{m=1}^\infty (m-1)t^{m-2} N_{[m]}(z)=
\left (\sum_{k=1}^\infty t^{k-1} JN_{[k]}(z)\right )
\left (\sum_{l=1}^\infty t^{l-1} N_{[l]}(z)\right )
\EQn
Comparing the coefficients of $t^{m-2}$ of the both sides 
of the equation above, we have
\BQn
(m-1) N_{[m]}(z)&=& \sum_{\substack{k+l=m\\ k, l\geq 1}}  
JN_{[k]}\cdot N_{[l]}
\EQn
for any $m\geq 2$. Hence we get Eq.\,(\ref{Nm}). 
\epfv

From Eq.\,(\ref{N1}), (\ref{Nm}) and by using the mathematical induction, 
it is easy to see that we have the following lemma.

\begin{lemma}\label{L3.4}
$(a)$ $o(N_{[m]})\geq m+1$ for any $m\geq 0$.

$(b)$ Suppose $H(z)\in \bC[z]^{\times n}$, then , for any $m\geq 1$, 
 $N_{[m]}\in \bC[z]^{\times n}$ with $\deg N_{[m]}(z) \leq (\deg H-1)m+1$.

$(c)$ If $H(z)$ is homogeneous  of degree $d$, then,
 $N_{[m]}(z)$ is homogeneous of degree
$(d-1)m+1$ for any $m\geq 1$.
\end{lemma}

Note that, by Lemma \ref{L3.4}, $(a)$, the infinite sum 
$\sum_{m=1}^\infty t_0^{m-1} N_{[m]}(z)$ makes sense for any
complex number $t=t_0$. In particular, when $t=1$, $G_{t=1}(z)$ gives us 
the formal inverse $G(z)$ of $F(z)$.

\begin{theo} \label{T3.5} {\bf (Recurrent Inversion Formula)} 

Let  $\{N_{[m]}(z) | m\geq 1 \}$ be the sequence defined 
by Eq. $(\ref{N1})$ and $(\ref{Nm})$ recurrently. Then   
the formal inverse of $F(z)=z-H(z)$ is given by 
\BQ\label{1-Inv}
G(z)=z+\sum_{m=1}^\infty N_{[m]}(z)
\EQ
\end{theo}

The following proposition also seems interesting to us. 
It basically says that $N_t(z)$ gives 
a family of formal maps from $\bC^n$ to $\bC^n$, 
which are ``closed'' under 
the inverse operation.

\begin{propo}\label{P3.10}
For any $s\in \bC$, 
 the formal inverse of $U_{s, t}(z)=z-sN_t(z)$
is given by
$V_{s, t}(z)=z+ s N_{t+s}(z)$.
Actually, $U_{s, t}= F_{t+s}\circ G_t$ and $V_{s, t}= F_{t}\circ G_{s+t}$.
\end{propo}
\pf 
\BQn
F_{t+s}\circ G_t &=& G_t(z)-(t+s)H(G_t)\\
&=& z+tN_t(z)-(t+s)N_t(z)\\
&=& z-sN_t(z)\\
&=&U_{s, t}(z)
\EQn
 Similarly, we can prove $V_{s, t}= F_{t}\circ G_{s+t}$.
\epfv

Another special property of $N_t(z)$ is given by the following proposition.

\begin{propo}
For any $U_{[0]}(z)\in \bC[[z]]^{\times n}$, the unique power series solution $U_t(z)$ 
in both $t$ and $z$ of the Cauchy problem
\BQ
\frac {\p U_t}{\p t} &=& JU_t N_t \label{GPDE}\\
U_{t=0}(z) &=& U_{[0]}(z)\label{GPDE-B}
\EQ  
is given by $U_t(z)=U_{[0]}(z+tN_t(z))$.
\end{propo}
\pf It is easy to see that the power series solution 
in both $t$ and $z$ of the Cauchy problem is unique. So it will be enough to 
check that $U_t(z)=U_{[0]}(z+tN_t(z))$ 
is such a solution.

\BQn
\frac {\p U_t}{\p t} &=&
\frac {\p}{\p t}U_{[0]}(z+tN_t(z))\\
&=& 
JU_{[0]}(z+tN_t(z))(N_t(z)+t\frac {\p N_t}{\p t}(z))\\
&=& JU_{[0]}(z+tN_t(z))(N_t(z)+tJN_t(z)N_t(z))\\
&=& JU_{[0]}(z+tN_t(z))(\text{I}+tJN_t(z)) N_t(z)\\
&=& J (U_{[0]}(z+tN_t(z))N_t(z)\\
&=& JU_t(z)N_t(z)
\EQn
\epfv

Next we derive more consequences of Theorem \ref{T3.2}.
\begin{propo}\label{New-P}
Let $F(z)=z-H(z)$ be a formal map with $o(H(z))\geq 2$. Then 
the inverse map $G(z)=z+H(z)$ if and only if $JH\cdot H=0$.
\end{propo}
\pf First, by Eq.\,(\ref{E3.1}) and Eq.\,(\ref{MainEQ}), we have
\begin{align*}
(JN_t\cdot N_t)(F_t)=(I-&tJH)^{-1}JH \cdot H.
\end{align*}
Hence, we have $JH \cdot H=0$ if and only if $JN_t\cdot N_t=0$.
    
Now we assume $JH\cdot H=0$. By Eq.\,(\ref{PDE}) and the fact above, 
we have $\frac {\p N_t}{\p t}=0$. By  Eq.\,(\ref{PDE-B}), 
we have $N_t(z)=H(z)$. 
Hence, the proposition follows.

Next we assume $G(z)=z+H(z)$ and to show $JH\cdot H=0$. 

First, from the equation,
\begin{align}
z=F(G(z))=z+H(z)-H(z+H(z)),
\end{align}
we see that 
\begin{align}\label{H=H}
H(z+H(z))=H(z). 
\end{align}

We consider the powers $F^{[m]}(z)$ $(m\geq 1)$ defined by
 Eq.\,(\ref{F-m>0}) of $F(z)=z-H(z)$ and
make the following claim.
\vskip3mm
\underline{\it Claim:} {\it For any $m\geq 1$, $F^{[m]}(z)=z-mH(z)$.}
\vskip3mm

\underline{\it Proof of Claim:} We use the mathematical induction on $m\geq 1$. 
The case $m=1$ is trivial. 
For $m>1$, we have
\begin{align*}
F^{[m]}(z)&=F^{[m-1]}(F(z))\\
&=z+H(z)-(m-1)H(z-H(z))\\
\intertext{Applying Eq.\,(\ref{H=H}):}
&=z-H(z)-(m-1)H(z)\\
&=z-mH(z).
\end{align*}

Hence the claim holds.

Next we consider the formal flow $F(z; t)$ of $F(z)$, i.e. the unique formal 
maps with coefficients in $\bC[t]$ such that
\begin{enumerate}
\item $F(z; 0)=z$ and $F(z; 1)=F(z)$.
\item For any $s, t\in \bC$, we have $F(F_t(z); s)=F(z; s+t)$.
\end{enumerate}

By Proposition $2.1$ and Lemma $2.4$ in \cite{Z1}, 
we know that there exists a unique $a(z)\in\bC[[z]]^{\times n}$ 
with $o(a(z))\geq 2$, 
which was called the D-Log of $F(z)$ in \cite{WZ}, such that
\begin{align}
e^{t a(z)\frac {\p}{\p z} } z=F(z, t).
\end{align}

Note that $F(z; m)=F^{[m]}(z)$ for any $m\geq 1$. Since $F(z; t)\in \bC[t][[z]]^{\times n}$, 
the coefficients of all monomials appearing in $F(z; t)$ are polynomials in $t$.
By the claim above, we see that the coefficients of all monomial appearing in $F(z; t)$ 
but not in $tH(z)$ vanish at any $m\geq 1$, hence they must be identically zero.
Therefore, we have
\begin{align}
F(z, t)=z-tH(z).\label{Short-Dlog}
\end{align}

From Eq.\,(\ref{Short-Dlog}) above, it is easy to see that  
\begin{align}
a(z)&=-H(z).\\
( a(z)\frac \p{\p z})^2 z &=H(z)\frac \p{\p z} H(z)=JH\cdot H=0.
\end{align}
Hence we are done.
\epfv

One immediate consequence of Proposition \ref{New-P} 
is the following theorem which was first proved in \cite{BCW}.

\begin{theo}\label{Thm-BCW} \cite{BCW} \,
Let $F(z)=z+H(z)$ be a polynomial map with $H(z)$ being homogeneous of degree $d\geq 2$. 
If $J(H)^2=0$, then the formal inverse map $G=z-H(z)$. 
\end{theo}
\pf When $H(z)$ is homogeneous of degree $d\geq 2$, by Euler's formula, we have
\begin{align}
JH^2(z)z=dJH\cdot H(z).
\end{align}
Hence, $JH^2(z)=0$ if and only $JH\cdot H=0$. Then the theorem follows immediately 
from Proposition \ref{New-P} above.
\epfv

In \cite{BCW}, H. Bass, E. Connell and D. Wright reduced the Jacobian 
conjecture to the cases when $H(z)$ is homogeneous of degree $3$. For further
reductions in this direction, see \cite{D1} and \cite{D3}. Note that, 
when $H(z)$ is homogeneous, 
the Jacobian condition $j(F)=1$ is equivalent to the condition that 
the Jacobian matrix $JH(z)$ is nilpotent. 
Next we will derive some consequences from our recurrent inversion formula 
for the case when $H(z)$ is homogeneous of degree $d$ ($d\geq 2$).

\begin{propo}\label{P3.6}
If $H$ is homogeneous of degree $d\geq 2$,  we have
\BQ
N_t(z)&=&\frac 1d JN_t(z)(z-(d-1)tN_t)\label{NtInJNt}\\
JN_t(z)\cdot z &=& d\left (I+\frac {(d-1)\, t}d JN_t(z)\right )N_t(z)
\label{JNz}
\EQ
\end{propo}

\pf First, by Euler's formula, we have $dH(z)=JH(z)z$. By composing
with $G_t$ from right, we get
\BQn
 dH(G_t)=JH(G_t) G_t(z)
\EQn
From Eq.\,(\ref{E3.2}), we have
\BQn
JH(G_t)=JN_t(z)J G_t^{-1}(z).
\EQn
Therefore, we have
\begin{align*}
dN_t(z)&= JN_t(z)J G_t^{-1}(z) (z+tN_t(z)) \\
dN_t(z)&= J G_t^{-1}(z) JN_t(z) (z+tN_t(z)) \\
\intertext{because $J G_t^{-1}(z)=(I+tJN_t(z))^{-1}$ commutes 
with the matrix $JN_t(z)$.}
d JG_t(z) N_t(z))&= JN_t(z) (z+tN_t(z))\\
d(1+tJN_t(z))N_t(z)&= JN_t(z)(z+tN_t(z))
\end{align*}

Solve $N_t(z)$ and $JN_t(z)z$ from the last equation above, 
we get Eq.\,(\ref{NtInJNt}) and (\ref{JNz}). 
\epfv

\begin{propo}\label{P3.7}
If $H$ is homogeneous of degree $d\geq 2$, then we have
\BQ
N_t(z)&=&\frac 1d \sum_{k=1}^\infty (-1)^{k-1}\left 
(\frac {(d-1)\, t}d\right )^{k-1} JN_t^k(z)\cdot z \label{E3.23}\\
\frac {\p N_t}{\p t}(z)&=&\frac 1d \sum_{k=1}^\infty (-1)^{k-1}\left 
(\frac {(d-1)\, t}d\right )^{k-1} JN_t^{k+1}(z)\cdot z\label{E3.24}
\EQ
\end{propo}

The reason that we think the proposition above 
is interesting is because that it 
writes $N_t(z)$ and $JN_t(z)$ in terms of $JN_t^k(z)\cdot z$ $(k\geq 1)$. 
In particular, when $JH(z)$ is nilpotent, $JN_t(z)$ is also nilpotent 
by Lemma \ref{L3.1} and the sums above are finite sums.

\pf First note that Eq.\,(\ref{E3.24}) can be obtained 
by multiplying $JN_t(z)$ from left to the both sides of Eq.\,(\ref{E3.23})
and then applying Eq.\,(\ref{PDE}). So we only need to show Eq.\,(\ref{E3.23}).
By Eq.\,(\ref{JNz}), we have
\BQn
N_t(z) = \frac 1d \left (I+\frac {(d-1)\, t}d JN_t(z)\right )^{-1}
JN_t(z)\cdot z
\EQn
But 
\BQn
\left (I+\frac {(d-1)\, t}d JN_t(z)\right )^{-1}=\sum_{k=1}^\infty (-1)^{k-1}\left 
(\frac {(d-1)\, t}d\right )^{k-1} JN_t^{k-1}(z)
\EQn
Hence, Eq.\,(\ref{E3.23}) follows immediately.
\epfv

Next, we give one more proof for Theorem  \ref{Thm-BCW}. 

\underline {\it Second Proof of Theorem \ref{Thm-BCW} }: 
First from Lemma \ref{L3.1}, $(b)$, we know that $JN_t^2(z)=0$.
By Eq.\,(\ref{E3.23}) and Euler's formula, we have
\BQn 
N_t(z)&=&\frac 1d JN_t(z)\cdot z \\
&=& \frac 1d \sum_{m=1}^\infty t^{m-1} JN_{[m]}(z) \cdot z\\
&=&\frac 1d \sum_{m=1}^\infty ( (d-1)m+1) t^{m-1} N_{[m]}(z)
\EQn
While $N_t(z)=\sum_{m=1}^\infty t^{m-1} N_{[m]}(z)$. Comparing the coefficients
of $t^{m-1}$, we get 
\BQn
N_{[m]}(z)=
\frac {(d-1)m+1}d N_{[m]}(z)\\
\frac {(d-1)(m-1)}d N_{[m]}(z)=0
\EQn
Hence $N_{[m]}(z)=0$ for any $m\geq 2$.
\epfv

Unfortunately, both proofs given in this section
for Theorem \ref{Thm-BCW}  
fail for the cases $JH^k(z)=0$ with $k\geq 3$.
At least from the proofs above, we can see that the cases 
$JH^k(z)=0$ $(k\geq 3)$ are dramatically different 
from the case $JH^2(z)=0$ and will
be much more difficult to study. D. Wright \cite{W1}
has shown that, for $n=3$ and $d=3$, 
the Jacobian conjecture is true.
Hubber \cite{H} proves the Jacobian conjecture  when 
$n=4$ and $d=3$.

Finally, we make the following remarks on the key partial 
differential equation (\ref{PDE}) involved in this section.

By Theorem \ref{T3.2} and Proposition \ref{P3.3}, 
it is east to see that, 
for a polynomial map $F(z)=z-H(z)$ with $H(z)$ homogeneous of 
degree $d$ ($d\geq 2$), $F(z)$ is an automorphism of $\bC^n$ 
if and only if the unique solution $N_t(z)$ of 
the Cauchy problem (\ref{PDE})-(\ref{PDE-B}) 
is a polynomial solution in both $t$ and $z$. 
Combining with 
the reduction theorem in \cite{BCW}, we see that the 
Jacobian conjecture is equivalent 
to the following conjecture.

\begin{conj}
Let $H(z)$ be homogeneous of 
degree $d=3$ with the Jacobian matrix $JH(z)$ nilpotent. Then the
unique solution $N_t(z)$ of the Cauchy problem
$(\ref{PDE})$-$(\ref{PDE-B})$ is a polynomial 
solution in both $t$ and $z$.
\end{conj}

Since the Jacobian conjecture has been proved  
by Wang \cite{Wa} for the case $d=2$, 
hence the statement in the conjecture above is also true 
in this case. 
It will be very interesting to
see a proof for this fact by using some PDE methods. 
We hope that some PDE approaches to the Cauchy problem 
(\ref{PDE})-(\ref{PDE-B}) will 
provide some new understanding to the Jacobian conjecture. 
It is interesting to notice the 
resemblance of the PDE (\ref{PDE}) with the well known 
inviscid Burgers equation in $n$ variables, 
which is the following
partial differential equation. 

\BQ\label{Burgers}
\frac {\p U_t}{\p t}(z) + (JU_t)^\tau U_t=0
\EQ
where $(JU_t)^\tau$ is 
the transpose of the Jacobian matrix of $U_t$ with respect to $z$.

Set $M_t(z)=-U_t(z)$, then Eq, (\ref{Burgers}) becomes
\BQ\label{Burgers1}
\frac {\p M_t}{\p t}(z) = (JM_t)^\tau M_t
\EQ
which differs from our equation (\ref{PDE}) only by  
the transpose for the Jacobian matrix of 
the unknown functions.  Interesting,  M. de Bondt, A. van den Essen
\cite{BE1} and G. Meng \cite{M} recently
have made a breakthrough 
on the Jacobian conjecture.
They reduced the Jacobian conjecture to 
polynomial maps $F(z)=z-H(z)$ with 
$H(z)=\nabla P(z)$ of some polynomials 
$P(z)\in\bC[z]$. Note that, in this case, 
$JH(z)=(\frac {\p^2 P(z)}{\p z_i\p z_j})$ 
is nothing but the Hessian matrix of $P(z)$. 
In particular, $JH(z)$ is symmetric and 
our equation (\ref{PDE}) does become exactly
the $n$-dimensional inviscid 
Burgers' equation Eq.\,(\ref{Burgers1}).
For further study for formal maps of the form 
$F(z)=z-H(z)$ with $H(z)=\nabla P(z)$ of some polynomials 
$P(z)\in\bC[z]$. 
See \cite{BE1}, \cite{BE2}, \cite{EW}, \cite{M}, 
\cite{Z2} and \cite{Z3}.

\renewcommand{\theequation}{\thesection.\arabic{equation}}
\renewcommand{\therema}{\thesection.\arabic{rema}}
\setcounter{equation}{0}
\setcounter{rema}{0}

\section{\bf A Differential-Free Recurrent Formula} \label{S4}

In this section, we derive a differential-free recurrent inversion formula 
for formal inverse $G(z)$ of the polynomial maps $F(z)=z-H(z)$ 
with $H(z)$ homogeneous of degree $d\geq 2$.
When $d=3$, our formula is same as the one given by 
Dru\.zkowski \cite{D2} except the symmetric 
multi-linear form involved is also given explicitly in our approach. 

In this section, we will always assume that $H(z)$ is homogeneous of 
degree $d$ for some $d\geq 2$. For any $1\leq j\leq d$, we set 
$U^{(j)}=(U^{(j)}_1, U^{(j)}_2, \cdots, U^{(j)}_n)$, 
 and
$D^{(j)}=\sum_{k=1}^n  U^{(j)}_k \frac {\p}{\p z_k}$,
where $U^{(j)}_i$ are formal variables.
For any $1\leq i\leq n$, we define 
the $d$-linear form $B_i$ by

\BQ\label{E4.1}
B_i( U^{(1)}, U^{(2)}, \cdots, U^{(d)})
=\frac 1{d!} D^{(1)}D^{(2)}\cdots D^{(d)} H_i(z).
\EQ

Note that the differential operators
 $D^{(j)}$ $(1\leq j\leq d)$ commute 
with each other, 
hence the $d$-linear forms 
$B_i( U^{(1)}, U^{(2)}, \cdots, U^{(d)})$ ($1\leq i\leq n$)
are symmetric. We set $B=(B_1, B_2, \cdots, B_n)$.

\begin{lemma} \label{L4.1}
Let $B( U^{(1)}, U^{(2)}, \cdots, U^{(d)})$ be the
symmetric multi-linear forms defined by Eq. $(\ref{E4.1})$. Then
\BQ\label{E4.2}
B(z, z, \cdots, z)=H(z).
\EQ
\end{lemma}
\pf
Note that, by Eq.\,(\ref{E4.1}) and Euler's formula, 
we have
\begin{align*}
& B_i(z, z, \cdots, z) \\
&= \frac 1{d!} \sum_{k_1, \cdots, k_d=1}^n z_{k_1}
z_{k_2}\cdots z_{k_d}\frac {\p^d H_i}{\p z_{k_1}
\p z_{k_2}\cdots \p z_{k_d}} \\
 &= \frac 1{d!}
\sum_{k_2, \cdots, k_d=1}^n 
z_{k_2}\cdots z_{k_d}\sum_{k_1=1}^n  z_{k_1} \frac {\p}{\p z_{k_1}} 
\left (\frac {\p^{d-1} H_i}{\p z_{k_2}
\p z_{k_3}\cdots \p z_{k_d}}\right ) \\
 &= \frac 1{d!}
\sum_{k_2, \cdots, k_d=1}^n 
z_{k_2}\cdots z_{k_d}
\frac {\p^{d-1} H_i}{\p z_{k_2}
\p z_{k_3}\cdots \p z_{k_d}} \\
\intertext{This is because that $\frac {\p^{d-1} H_i}{\p z_{k_2}
\p z_{k_3}\cdots \p z_{k_d}}$ is homogeneous of degree $1$.}
&= \frac 1{d!}
\sum_{k_3, \cdots, k_d=1}^n 
z_{k_3}\cdots z_{k_d}\sum_{k_2=1}^n  z_{k_2} \frac {\p}{\p z_{k_2}} 
\left (\frac {\p^{d-2} H_i}{\p z_{k_3}
\p z_{k_3}\cdots \p z_{k_d}}\right ) \\
 &= \frac 2{d!}
\sum_{k_3, \cdots, k_d=1}^n 
z_{k_3}\cdots z_{k_d}
\frac {\p^{d-2} H_i}{\p z_{k_3}
\p z_{k_3}\cdots \p z_{k_d}} 
\intertext{This is because that $\frac {\p^{d-2} H_i}{\p z_{k_3}
\p z_{k_3}\cdots \p z_{k_d}}$ is homogeneous of degree $2$.}
&=\cdots
\end{align*}
By repeating the procedure above, it is easy to see that Eq.\,(\ref{E4.2})
holds.
\epfv

By Lemma \ref{L3.4} $(c)$, we can 
write the formal inverse $G(z)$ as 
$G(z)=z+\sum_{m=1}^\infty N_{[m]}(z)$, where 
$N_{[m]}(z)$ $(m\geq 1)$ are
homogeneous of degree $m(d-1)+1$. Note that, 
for convenience, we have set $N_{[0]}(z)=z$.
Now we can write down our second recurrent inversion formula as follows.

\begin{theo}\label{T4.2}
We have the following recurrent formula for the formal inverse 
$G(z)=z+\sum_{m=1}^\infty N_{[m]}(z)$.
\BQ
N_{[1]}(z)&=& H(z)\\
N_{[m+1]} (z)&=&\sum_{\substack {k_1+\cdots +k_n=m\\k_i\geq 0}} 
B(N_{[k_1]}, N_{[k_2]}, \cdots, N_{[k_n]}) \label{2-Inv}
\EQ 
for any $m\geq 1$.
\end{theo}
\pf By replacing $z$ with $G(z)$ in Eq.\,(\ref{E4.2}), 
we get
\BQn
B(\sum_{m=0}^\infty N_{[m]}, \sum_{m=0}^\infty N_{[m]}, \cdots, 
\sum_{m=0}^\infty N_{[m]})
=H(G(z))
\EQn
While, by Eq.\,(\ref{E3.2}) in Lemma \ref{L3.1}, we have
\BQn
H(G(z))=N(z)=\sum_{m=1}^\infty N_{[m]}(z)
\EQn
By using the equations above and the fact that 
$B$ is a multi-linear form, we get
\BQn
\sum_{m=1}^\infty N_{[m]}(z)&=&
B(\sum_{m=0}^\infty N_{[m]}, 
\sum_{m=0}^\infty N_{[m]}, \cdots, \sum_{m=0}^\infty N_{[m]})\\
&=&
\sum_{r= 0}^\infty 
\sum_{\substack {k_1+\cdots +k_d=r\\k_i\geq 0}} 
B(N_{[k_1]}, N_{[k_2]}, \cdots, N_{[k_d]})
\EQn
By comparing the homogeneous parts of both sides of the equation above,
we get Eq.\,(\ref{2-Inv}).
\epfv

{\small \sc Department of Mathematics, Illinois State University, Normal IL 61790-4520.}

{\em E-mail}: 
wzhao@ilstu.edu.

\end{document}